\documentclass[a4paper,12pt]{article}
\usepackage{tipa,amsmath,amsfonts,amssymb,amsthm,amscd,mathrsfs,latexsym,bm}
\usepackage[all]{xy}
\usepackage{graphicx}
\usepackage{mathptm,pslatex}

 \newcommand{\lam}{\lambda}
 
\newcommand{\ft}{\mathfrak t}

\newcommand{\fs}{\mathfrak s} 
 
\DeclareMathOperator{\Std}{Std}

\def\qed{\hfill$\Box$\smallskip}

\newtheorem{prop}{Proposition}[section]
\newtheorem{thm}[prop]{Theorem}\newtheorem{cor}[prop]{Corollary}
\newtheorem{lem}[prop]{Lemma}\newtheorem{dfn}[prop]{Definition}

\title{Jucys-Murphy elements and centers of cellular algebras
\thanks {keywords: Jucys-Murphy elements, cellular algebras, center.}}
\author{Yanbo Li\\[5pt]Department of Information and Computing Sciences,
\\Northeastern University at Qinhuangdao; \\Qinhuangdao, 066004, P.R. China \\[5pt]School of Mathematics Sciences,
Beijing Normal University;\\
Beijing, 100875, P.R. China\\[1.5pt] E-mail: liyanbo707@163.com}
\begin{document}
\maketitle
\begin{abstract}

Let $R$ be an integral domain and $A$ a cellular
algebra over $R$ with a cellular basis $\{C_{S,T}^\lam \mid \lam\in\Lambda, S,T\in
M(\lam)\}$. Suppose that $A$ is equipped with a family
of Jucys-Murphy elements which satisfy the separation condition in the sense of A.
Mathas \cite{M}. Let $K$ be the field of fractions of $R$ and
$A_{K}=A\bigotimes_{R}K$. We give a necessary and sufficient
condition under which the center of $A_{K}$ consists of the symmetric polynomials in Jucys-Murphy
elements.
\end{abstract}

\section{Introduction}

Jucys-Murphy elements were constructed for the group algebras of
symmetric groups first. The combinatorics of these elements allow
one to compute simple representations explicitly and often easily in
the semisimple case. Then Dipper, James and Murphy \cite{DG1},
\cite{DG2}, \cite{DG3}, \cite{DG4}, \cite{DGM} did a lot of work on
representations of Hecke algebras and produced analogues of
the Jucys-Murphy elements for Hecke algebras of types A and
B. The constructions for other algebras can be found in \cite{HR},
\cite{RA2} and so on. In \cite{DG2}, Dipper and James conjectured
that the center of a Hecke algebra of type A consists of symmetric
polynomials in Jucys-Murphy elements. The conjecture was proved by Francis and Graham
\cite{FG} in 2006. In \cite{B}, Brundan proved that the center of
each degenerate cyclotomic Hecke algebra consists of symmetric
polynomials in the Jucys-Murphy elements. An analogous conjecture
for Ariki-Koike Hecke algebra is open in non-semisimple case.

Cellular algebras were introduced by Graham and Lehrer \cite{GL} in
1996, motivated by previous work of Kazhdan and Lusztig \cite{KL}.
The theory of cellular algebras provides a systematic
framework for studying the representation theory of non-semisimple
algebras which are deformations of semisimple ones. Many classes of algebras from
mathematics and physics are found to be cellular, including Hecke
algebras of finite type, Ariki-Koike Hecke algebras, $q$-Schur algebras,
Brauer algebras, partition algebras, Birman-Wenzl algebras
and so on, see \cite{G}, \cite{GL}, \cite{Xi1},
\cite{Xi2} for details.

The fact that most of the algebras which have Jucys-Murphy elements
are cellular leads one to defining Jucys-Murphy elements for general
cellular algebras. In \cite{M}, Mathas did some work in this
direction. By the definition of Mathas, we investigate the
relations between the centers and the Jucys-Murphy elements of
cellular algebras.

Let $R$ be an integral domain and $A$ a cellular $R$-algebra with a
cellular basis $\{C_{S,T}^\lam \mid S,T\in
M(\lam),\lam\in\Lambda\}$. Let $K$ be the field of fractions of $R$
and $A_{K}=A\bigotimes_{R}K$.  Suppose that $A$ is equipped with a
family of Jucys-Murphy elements $L_{1},\ldots,L_{m}$ which satisfy
the separation condition \cite{M}. For any
$\lam\in\Lambda$, $\{c_{\lam}(i)\mid 1\leq i\leq m\}$ is a family of
elements in $R$. Then the main result of this paper is the following theorem.

\bigskip

\noindent{\bf Theorem.} \emph{Suppose that every symmetric
polynomial in $L_{1},\ldots,L_{m}$ belongs to the center of $A_{K}$.
Then the following are equivalent.}

\noindent(1) \emph{The center of $A_{K}$ consists of symmetric
polynomials in Jucys-Murphy elements.}

\noindent(2) \emph{$\{c_{\lam}(i)\mid 1\leq i\leq m\}$ can not be
obtained from $\{c_{\mu}(i)\mid 1\leq i\leq m\}$ by permutations for
arbitrary $\lam,\mu\in\Lambda$ with $\lam\neq\mu$.}

\bigskip

The condition in the above theorem is also a necessary condition
for the center of $A$ consisting of the symmetric polynomials in
Jucys-Murphy elements. Moreover, by this theorem, we can prove that
the centers of Ariki-Koike Hecke algebras consist of the symmetric
polynomials in Jucys-Murphy elements in semisimple case. The proof
is different from Ariki's in \cite{A} and A. Ram's in \cite{RA}.

\section{Cellular algebras and Jucys-Murphy elements}
In this section, we first recall the definition of cellular algebras
 and then give a quick review of the results under the
so-called separation condition in A. Mathas' paper \cite{M}.

\begin{dfn} {\rm (\cite{GL} 1.1)}
\label{dfgl1}Let $R$ be a commutative ring with identity. An
associative unital $R$-algebra is called a cellular algebra with
cell datum $(\Lambda, M, C, i)$ if the following conditions are
satisfied:

(C1) The finite set $\Lambda$ is a poset. Associated with each
$\lam\in\Lambda$, there is a finite set $M(\lam)$. The algebra $A$
has an $R$-basis $\{C_{S,T}^\lam \mid S,T\in
M(\lam),\lam\in\Lambda\}$.

(C2) The map $i$ is an $R$-linear anti-automorphism of $A$ with
$i^{2}=id$ which sends $C_{S,T}^\lam$ to $C_{T,S}^\lam$.

(C3) If $\lam\in\Lambda$ and $S,T\in M(\lam)$, then for any element
$a\in A$, we have\\
$$aC_{S,T}^\lam\equiv\sum_{S^{'}\in
M(\lam)}r_{a}(S^{'},S)C_{S^{'},T}^{\lam} \,\,\,\,(\rm {mod}\,\,\,
 A(<\lam)),$$ where $r_{a}(S^{'},S)\in R$ is independent of $T$ and
where $A(<\lam)$ is the $R$-submodule of $A$ generated by
$\{C_{S^{''},T^{''}}^\mu \mid S^{''},T^{''}\in M(\mu),\mu<\lam\}$.

Apply $i$ to the equation in (C3), we obtain

$(C3^{'})\,\, C_{T,S}^\lam i(a)\equiv\sum\limits_{S^{'}\in
M(\lam)}r_{a}(S^{'},S)C_{T,S^{'}}^{\lam} \,\,\,\,(\rm mod
\,\,\,A(<\lam)).$
\end{dfn}

Let $R$ be an integral domain. Given a cellular algebra $A$, we will
also assume that $M(\lam)$ is a poset with an order $\leq_{\lam}$.
Let $ M(\Lambda)=\bigsqcup_{\lam\in \Lambda}M(\lam)$, we consider
$M(\Lambda)$ as a poset with an order $\leq$ as follows.

$$
S\leq T\Leftrightarrow\begin{cases} S\leq_{\lam} T,&\text{if $S,T\in M(\lam)$;}\\
\lam\leq\mu,\, &\text{if $S\in M(\lam)$, $T\in M(\mu)$.}\end{cases}
$$

Let $K$ be the field of fractions of $R$ and
$A_{K}=A\bigotimes_{R}K$. We will consider $A$ as a subalgebra of
$A_{K}$.

\begin{dfn} {\rm (\cite{M} 2.4)}\label{dfjm} Let $R$ be an integral domain and $A$ a cellular
algebra. A family of elements $L_{1},\ldots,L_{m}$ are called
Jucys-Murphy elements of $A$ if\\
{\rm (1)} $L_{i}L_{j}=L_{j}L_{i}$, for $1\leq i,j\leq m$;\\
{\rm (2)} $i(L_{j})=L_{j}$, for $j=1,\cdots,m$;\\
{\rm (3)} For all $\lam\in\Lambda$, $S,T\in M(\lam)$ and $L_{i}$, $i=
1,\cdots, m$, $$C_{S,T}^{\lam}L_{i}\equiv
c_{T}(i)C_{S,T}^{\lam}+\sum_{V<
T}r_{L_{i}}(V,T)C_{S,V}^{\lam}\,\,\,\,(\mod\,\,A({<\lam})),$$ where
$c_{T}(i)\in R, \,\,r_{L_{i}}(T,V)\in R$. We call $c_{T}(i)$ the
content of $T$ at $i$. Denote $\{c_{T}(i)\mid T\in M(\Lambda)\}$ by $\mathscr{C}(i)$ for $i=1, 2, \cdots, m$.
\end{dfn}

\noindent{\bf Example.}\, Let $K$ be a field. Let $A$ be the group
algebra of symmetric group $S_n$. Set
$L_{i}=\sum\limits_{j=1}^{i-1}(i,j)$ for $i=2,\cdots,n$. Then
$L_{i},\,\,\, i=2,\cdots,n$, is a family of Jucys-Murphy elements of $A$.

\begin{dfn} {\rm (\cite{M} 2.8)}
Let $A$ be a cellular algebra with Jucys-Murphy elements
$\{L_{1},\ldots,L_{m}\}$. We say that the Jucys-Murphy elements
satisfy the separation condition if for any $S,T\in M(\Lambda),
S\leq T, S\neq T$, there exists some $i$ with $1\leq i\leq m$ such that
$c_{S}(i)\neq c_{T}(i)$.
\end{dfn}
\noindent{\bf Remark.} The separation condition forces $A_{K}$ to be
semisimple (c.f. \cite {M}).

From now on, we always assume that $A$ is a cellular algebra
equipped with a family of Jucys-Murphy elements which satisfy the
separation condition. We now recall some results of \cite{M}.

\begin{dfn} {\rm (\cite{M} 3.1)} Let $A$ be a cellular algebra with Jucys-Murphy elements
$\{L_{1},\ldots,L_{m}\}$. For $\lam\in\Lambda$, $S,T\in M(\lam)$,
define
$$F_{T}=\prod_{i}\prod_{c\in \mathscr{C}(i),c\neq c_{T}(i)}(L_{i}-c)/(c_{T}(i)-c)$$
and $f_{S,T}^{\lam}=F_{S}C_{S,T}^{\lam}F_{T}$.
\end{dfn}
Note that the coefficient of $C_{S,T}^{\lam}$ in the expansion of
$f_{S,T}^{\lam}$ is $1$ for any $\lam\in\Lambda$ and $S,T\in
M(\lam)$, see \cite{M} 3.3 (a). Then Mathas proved the following theorems.
\begin{thm} {\rm (\cite{M} 3.7)} Let $A$ be a cellular algebra with Jucys-Murphy elements
$\{L_{1},\ldots,L_{m}\}$. Let $\lam, \mu \in \Lambda$, $S,T\in
M(\lam)$ and $U,V\in M(\mu)$. Then\\ {\rm (1)}$$
f_{S,T}^{\lam}f_{U,V}^{\mu}=\begin{cases} \gamma_{T}f_{S,V}^{\lam},&\text{$\lam=\mu,\,\,\,T=U$,}\\
0,\, &\text{otherwise,}\end{cases}
$$
where $\gamma_{T}\in K$ and $\gamma_{T}\neq 0$ for all $T\in
M(\Lambda)$.\\
{\rm (2)} $\{f_{S,T}^{\lam}\mid S,T\in M(\lam), \lam \in \Lambda\}$ is a
cellular basis of $A_{K}$.\qed
\end{thm}

\begin{thm}\label{thm} {\rm (\cite{M} 3.16)} Let $A$ be a cellular algebra with Jucys-Murphy elements
$\{L_{1},\ldots,L_{m}\}$. Then\\
{\rm (1)} Let $\lam \in \Lambda$ and $T\in M(\lam)$. Then $F_{T}$ is a
primitive idempotent in $A_{K}$. Moreover, $\{F_{T}\mid T\in
M(\lam)\}$ is a complete set of
pairwise orthogonal primitive idempotents in $A_{K}$.\\
{\rm (2)} $F_{\lam}=\sum\limits_{T\in M(\lam)}F_{T}$ is a central
idempotent in $A_{K}$ for any $\lam\in\Lambda$. Moreover,
$\{F_{\lam}\mid \lam\in \Lambda\}$ is a complete set of
central idempotents which are primitive in $Z(A_{K})$.\\
{\rm (3)} In particular,
$1=\sum\limits_{\lam\in\Lambda}F_{\lam}=\sum\limits_{T\in
M(\Lambda)}F_{T}$ and $L_{i}=\sum\limits_{T\in
M(\Lambda)}c_{T}(i)F_{T}$.\qed
\end{thm}

\smallskip\bigskip

\section{Jucys-Murphy elements and centers of cellular algebras}
In \cite{M}, A. Mathas gave a relation between the center and
Jucys-Murphy elements of a cellular algebra.
\begin{prop} {\rm (\cite{M} 4.13)}\label{propjm} Let $A$ be a cellular algebra with Jucys-Murphy elements
$\{L_{1},\ldots,L_{m}\}$. For any $\lam\in\Lambda$ and $S,T\in
M(\lam)$, if $\{c_{S}(i)\mid 1\leq i\leq m\}$ can be obtained by
permutations from $\{c_{T}(i)\mid 1\leq i\leq m\}$, then every
symmetric polynomial in $L_{1},\ldots,L_{m}$ belongs to the center
of $A_{K}$.\qed
\end{prop}

In fact, the inverse proposition also holds.
\begin{prop}\label{propjm2} Let $A$ be a cellular algebra with a family of Jucys-Murphy elements
$\{L_{1},\ldots,L_{m}\}$. Suppose that every symmetric polynomial in
$L_{1},\ldots,L_{m}$ belongs to the center of $A_{K}$. Let
$\lam\in\Lambda$ and $S,T\in M(\lam)$. Then $\{c_{S}(i)\mid 1\leq
i\leq m\}$ can be obtained by permutations from $\{c_{T}(i)\mid
1\leq i\leq m\}$.
\end{prop}
\noindent {Proof:}\, Suppose that there exists some $\lam\in\Lambda$
and $S,T\in M(\lam)$ such that $\{c_{S}(i)\mid 1\leq i\leq m\}$ can
not be obtained by permutations from $\{c_{T}(i)\mid 1\leq i\leq
m\}$. Then there exists a symmetric polynomial $p$ such that
$$p(c_{S}(1),\ldots,c_{S}(m))\neq p(c_{T}(1),\ldots,c_{T}(m)).$$
Note that $L_{i}=\sum\limits_{X\in M(\Lambda)}c_{X}(i)F_{X}$, then
$$p(L_{1}, \cdots, L_{m})=\sum_{U\in M(\Lambda)}p(c_{U}(1),\ldots,c_{U}(m))F_{U}.$$
Multiply by $F_{T}$ on both sides, we get $p(L_{1}, \cdots, L_{m})F_{T}= p(c_{T}(1),\ldots,c_{T}(m))F_{T}$ from Theorem \ref{thm} (1), the equation $p(L_{1}, \cdots, L_{m})F_{S}= p(c_{S}(1),\ldots,c_{S}(m))F_{S}$ is obtained similarly.

On the other hand, since $p(L_{1}, \cdots, L_{m})\in Z(A_{K})$, then by Theorem \ref{thm} (3), $p(L_{1}, \cdots, L_{m})=\sum\limits_{\lam\in\Lambda}r_{\lam}F_{\lam}$, where $r_{\lam}\in K$.
Multiply by $F_{T}$ on both sides, we get $p(L_{1}, \cdots, L_{m})F_{T}=r_{\lam}F_{T}$. The equation $p(L_{1}, \cdots, L_{m})F_{S}=r_{\lam}F_{S}$ can be obtained similarly. Then $p(c_{T}(1),\ldots,c_{T}(m))=r_{\lam}=p(c_{S}(1),\ldots,c_{S}(m))$. It is a contradiction.\qed

By the above proposition, if every symmetric polynomial in
$L_{1},\ldots,L_{m}$ belongs to the center of $A_{K}$, then for any
$\lam\in\Lambda$ and $S,T\in M(\lam)$, we have $\{c_{S}(i)\mid 1\leq
i\leq m\}$ and $\{c_{T}(i)\mid 1\leq i\leq m\}$ are the same if we
do not consider the order. So we can denote any of them by
$\{c_{\lam}(i)\mid 1\leq i\leq m\}$.

Now we are in a position to give the main result of this paper.

\begin{thm}Let $R$ be an integral
domain and $A$ a cellular $R$-algebra with a cellular basis
$\{C_{S,T}^\lam \mid S,T\in M(\lam),\lam\in\Lambda\}$. Let $K$ be
the field of fractions of $R$ and $A_{K}=A\bigotimes_{R}K$. Suppose
that $A$ is equipped with a family of Jucys-Murphy elements
$L_{1},\ldots,L_{m}$ which satisfy the separation condition and all
symmetric polynomials in $L_{1},\ldots,L_{m}$ belong to the center
of $A_{K}$. Then the following are equivalent.\\
(1) The center of $A_{K}$ consists of all symmetric polynomials in the
Jucys-Murphy elements.\\
(2) For any $\lam,\mu\in\Lambda$ with $\lam\neq\mu$,
$\{c_{\lam}(i)\mid 1\leq i\leq m\}$ can not be obtained from
$\{c_{\mu}(i)\mid 1\leq i\leq m\}$ by permutations.
\end{thm}

To prove this theorem, we need the following two lemmas.
\begin{lem}\label{lmjm1} Let $X_{1}, X_{2}, \cdots, X_{m}$ be indeterminates over a field $K$ and
let $\{x_{1},\ldots,x_{m}\}$ and $\{y_{1},\ldots,y_{m}\}$ be two
families of elements in $K$. Suppose that there exists some $k\in
K$, such that $p(x_{1},\ldots,x_{m})=kp(y_{1},\ldots,y_{m})$ for any symmetric polynomial $p(X_{1}, X_{2}, \cdots, X_{m})\in K[X_{1}, X_{2}, \cdots, X_{m}]$. Then
$\{x_{1},\ldots,x_{m}\}$ can be obtained by permutations from
$\{y_{1},\ldots,y_{m}\}$.
\end{lem}
\noindent {Proof:}\, Clearly, if $p$ is a symmetric polynomial, then
$p^{2}$ is also a symmetric polynomial. Then
$$(p(x_{1},\ldots,x_{m}))^{2}=(kp(y_{1},\ldots,y_{m}))^{2}=k(p(y_{1},\ldots,y_{m}))^{2}.$$
Hence $(k^{2}-k)(p(y_{1},\ldots,y_{m}))^{2}=0$. Then $k^{2}-k=0$
since $p$ is arbitrary. So we have $k=0$ or $k=1$. If $k=0$, then
$p(x_{1},\ldots,x_{m})=0$ for any $p$. This is impossible. Then
$k=1$. That is $p(x_{1},\ldots,x_{m})=p(y_{1},\ldots,y_{m})$ for
arbitrary $p$.\qed

Let $\{k_{11},\ldots,k_{1m}\}$, $\ldots$, $\{k_{n1},\ldots,k_{nm}\}$
be $n$ families of elements in $K$ and $p$ a symmetric polynomial.
We will denote $p(k_{i1},\ldots,k_{im})$ by $p(i)$.
\begin{lem}\label{lmjm2} Suppose that  $\{k_{11},\ldots,k_{1m}\}$, $\ldots$, $\{k_{n1},\ldots,k_{nm}\}$ are
 $n$ families of elements in a field $K$ and $X_{1}, X_{2}, \cdots, X_{m}$ indeterminates.
Let $p_{1}^{'}(X_{1}, X_{2}, \cdots, X_{m})$, $\ldots$, $p_{n}^{'}(X_{1}, X_{2}, \cdots, X_{m})$ $\in K[X_{1}, X_{2}, \cdots, X_{m}]$ be $n$ symmetric polynomials such that
$$\begin{vmatrix} p_{1}^{'}(1) & \ldots &
p_{1}^{'}(n)\\\ldots & \ldots & \ldots\\p_{n}^{'}(1) & \ldots &
p_{n}^{'}(n)
\end{vmatrix}\neq 0.$$ Then there exist $n$ symmetric polynomials
$p_{1},\ldots,p_{n}$ such that
$$\begin{vmatrix} p_{1}(1) & p_{1}(2) & \ldots & p_{1}(n)\\0 & p_{2}(2) & \ldots & p_{2}(n)\\
\ldots & \ldots & \ldots &\ldots\\0 & 0 & \ldots &
p_{n}(n)\end{vmatrix}\neq 0.$$
\end{lem}
\noindent {Proof:}\, Without loss of generality, we assume that
$p_{1}^{'}(1)\neq 0$ and set $p_{1}=p_{1}^{'}$. Then let
$p_{2}=p_{2}^{'}-\dfrac{p_{2}^{'}(1)}{p_{1}(1)}p_{1}$. Clearly,
$p_{2}$ is a symmetric polynomial and $p_{2}(1)=0$. Moreover,
$$\begin{vmatrix} p_{1}(1) & p_{1}(2) & \ldots & p_{1}(n)\\0 & p_{2}(2) & \ldots & p_{2}(n)\\
p_{3}^{'}(1) & p_{3}^{'}(2) & \ldots & p_{3}^{'}(n)\\\ldots & \ldots & \ldots &\ldots\\
p_{n}^{'}(1) & p_{n}^{'}(2) & \ldots &
p_{n}^{'}(n)\end{vmatrix}=\begin{vmatrix} p_{1}^{'}(1) & \ldots &
p_{1}^{'}(n)\\\ldots & \ldots & \ldots\\p_{n}^{'}(1) & \ldots &
p_{n}^{'}(n)
\end{vmatrix}.$$
Repeat the above process similarly, we can find $p_{1},\ldots,p_{n}$
such that
$$\begin{vmatrix} p_{1}(1) & p_{1}(2) & \ldots & p_{1}(n)\\0 & p_{2}(2) & \ldots & p_{2}(n)\\
\ldots & \ldots & \ldots &\ldots\\0 & 0 & \ldots &
p_{n}(n)\end{vmatrix}=\begin{vmatrix} p_{1}^{'}(1) & \ldots &
p_{1}^{'}(n)\\\ldots & \ldots & \ldots\\p_{n}^{'}(1) & \ldots &
p_{n}^{'}(n)
\end{vmatrix}\neq 0.$$\qed

\noindent{\bf Proof of Theorem}. Since $p(L_{1},L_{2},\cdots,L_{m})=\sum\limits_{\lam\in\Lambda}p(c_{\lam}(1),\cdots,c_{\lam}(m))F_{\lam}$(see the proof of Proposition 4.12 in \cite{M}), then $(1)\Rightarrow (2)$ is obvious.
Now we prove $(2)\Rightarrow (1)$ by induction on the number of the
elements in the poset $\Lambda$. Denote the number by
$\sharp(\Lambda)$ and denote the elements in $\Lambda$ by natural numbers.

It is easy to know that we only need to find symmetric polynomials $p^{'}_{1}, p^{'}_{2}, \cdots, p^{'}_{n}$ such that
$$\begin{vmatrix} p^{'}_{1}(1) & p^{'}_{1}(2) & \ldots & p^{'}_{1}(n)\\p^{'}_{2}(1) & p^{'}_{2}(2) & \ldots & p^{'}_{2}(n)\\
\ldots & \ldots & \ldots &\ldots\\p^{'}_{n}(1) & p^{'}_{n}(2) & \ldots &
p^{'}_{n}(n)\end{vmatrix}\neq 0$$ where $n=\sharp(\Lambda)$.

For $\sharp(\Lambda)=1$, it is clear.

We now assume that $(2)\Rightarrow (1)$ holds for
$\sharp(\Lambda)=n$. Then by Lemma \ref{lmjm2}, there exist
symmetric polynomials $p_{1},\ldots,p_{n}$ such that
$$\begin{vmatrix} p_{1}(1) & p_{1}(2) & \ldots & p_{1}(n)\\0 & p_{2}(2) & \ldots & p_{2}(n)\\
\ldots & \ldots & \ldots &\ldots\\0 & 0 & \ldots &
p_{n}(n)\end{vmatrix}\neq 0.$$ We now assume that for any symmetric
polynomial $p$,
$$d:=\begin{vmatrix} p_{1}(1) & p_{1}(2) & \ldots & p_{1}(n) & p_{1}(n+1)\\0 & p_{2}(2) & \ldots & p_{2}(n) & p_{2}(n+1)\\
\ldots & \ldots & \ldots & \ldots & \ldots\\0 & 0 & \ldots &
p_{n}(n) & p_{n}(n+1)\\p(1) & p(2) & \ldots & p(n) &
p(n+1)\end{vmatrix}=0.$$ Then
$p(n+1)=k_{1}p(1)+k_{2}p(2)+\ldots+k_{n}p(n)$, where $k_{i}\in K$ is independent of $p$
for $i=1,\ldots,n$. Then we have
$p_{n}p(n+1)=k_{1}p_{n}p(1)+k_{2}p_{n}p(2)+\ldots+k_{n}p_{n}p(n)$ since
$p_{n}p$ is also a symmetric polynomial. Assume that $p_{n}(n+1)\neq
0$, then $p_{n}(n+1)p(n+1)=k_{n}p_{n}(n)p(n)$, or $p(n+1)=kp(n)$,
where $k\in K$ is independent of the choice of $p$. This implies that $p(n+1)=p(n)$ by Lemma
\ref{lmjm1}. It is a contradiction. Then $p_{n}(n+1)=0$. That
is $k_{n}p_{n}(n)p(n)=0$. Since $p_{n}(n)\neq 0$ and $p$ is arbitrary, then $k_{n}=0$. Repeat this process similarly, we have
$k_{i}=0$ for $i=1,\cdots,n$ and then $p(n+1)=0$. It is impossible
for $p$ is arbitrary. Then there exists a symmetric polynomial $p$
such that $d\neq 0$. This completes the proof.\qed

\begin{cor} Let $R$ be an integral domain and $A$ a cellular
algebra. Suppose that $A$ is equipped with a family
of Jucys-Murphy elements which separate $A$. If the center of $A$ consists of symmetric polynomials in Jucys-Murphy elements, then $\{c_{\lam}(i)\mid 1\leq i\leq m\}$ can not be
obtained from $\{c_{\mu}(i)\mid 1\leq i\leq m\}$ by permutations for
arbitrary $\lam,\mu\in\Lambda$ with $\lam\neq\mu$.\qed
\end{cor}

\smallskip\bigskip
\section{An application on Ariki-Koike Hecke algebras}

In this section, we prove that the center of a semisimple
Ariki-Koike Hecke algebra ($q\neq 1$) consists of the symmetric
polynomials in Jucys-Murphy elements. It is a new proof different
from Ariki's in \cite{A} and A. Ram's \cite{RA}.

Firstly, we recall some notions of combinatorics. Recall that a
partition of $n$ is a non-increasing sequence of non-negative integers $\lam=(\lam_1,\cdots,\lam_r)$ such that
$\sum_{i=1}^{r}\lam_i=n$. The diagram of a partition $\lam$ is the
subset $[\lam]=\{(i,j)\mid 1\leq j\leq\lam_{i}, i\geq 1\}$. The
elements of $\lam$ are called nodes. Define the residue of the node $(i,j)\in[\lam]$ to be $j-i$.
For any partition $\lam=(\lam_1,\lam_2,\cdots)$, the conjugate of
$\lam$ is defined to be a partition
$\lam'=(\lam'_1,\lam'_2,\cdots)$, where $\lam'_j$ is equal to the
number of nodes in column $j$ of $[\lam]$ for $j=1,2,\cdots$. For
partitions, we have the following simple lemma.

\begin{lem}\label{lmp}
Let $\lam$ and $\mu$ be two partitions of $n$. Then $\lam=\mu$ if
and only if all residues of nodes in $[\lam]$ and $[\mu]$ are the
same.\qed
\end{lem}

Given two partitions $\lam$ and $\mu$ of $n$, write
$\lam\trianglerighteq\mu$ if
$$\sum_{i=1}^{j}\lam_{i}\geq\sum_{i=1}^{j}\mu_{i},\,\,\,\, \text{for all $i\geq
1$}.$$ This is the so-called dominance order. It is a partial order.

A $\lam$-tableau is a bijection $\ft:
[\lam]\rightarrow\{1,2,\cdots,n\}$. We say $\ft$ a standard
$\lam$-tableau if the entries in $\ft$ increase from left to right
in each row and from top to bottom in each column. Denote by
$\ft^{\lam}$ (resp., $\ft_{\lam}$) the standard $\lam$-tableau, in
which the numbers $1,2,\cdots,n$ appear in order along successive
rows (resp., columns), The row stabilizer of $\ft^{\lam}$, denoted
by $S_{\lam}$, is the standard Young subgroup of $S_n$ corresponding
to $\lam$. Let $\Std(\lam)$ be the set of all standard
$\lam$-tableaux.

For a fixed positive integer $m$, a $m$-multipartitions of $n$ is an
$m$-tuple of partitions which sum to $n$. Let
$$\lam=((\lam_{11},\lam_{12},\cdots,\lam_{1i_{1}}),(\lam_{21},\lam_{22},\cdots,
\lam_{2i_{2}}),\cdots,(\lam_{m1},\lam_{m2},\cdots,\lam_{mi_{m}}))$$
be a $m$-multipartitions of $n$, we denote
$\lam_{j1}+\lam_{j2}+\cdots+\lam_{ji_{j}}$ by $n_{j\lam}$ for $1\leq
j\leq m$. A standard $\lam$-tableau is an $m$-tuple of standard
tableaux. We can define $\ft^{\lam}$ similarly.

\smallskip

Let $R$ be an integral domain, $q, u_{1}, u_{2}, \cdots, u_{m}\in R$
and $q$ invertible. Fix two positive integers $n$ and $m$. Then
Ariki-Koike algebra $\mathscr{H}_{n,m}$ is the associative
$R$-algebra with generators $T_{0}, T_{1}, \cdots, T_{n-1}$ and
relations
$$\begin{aligned}
&(T_{0}-u_{1})(T_{0}-u_{2})\cdots (T_{0}-u_{m})=0,\\
&T_{0}T_{1}T_{0}T_{1}=T_{1}T_{0}T_{1}T_{0},\\
&(T_i-q)(T_i+1)=0,\quad\text{for $1\leq i\leq n-1$,}\\
&T_iT_{i+1}T_i=T_{i+1}T_{i}T_{i+1},\quad\text{for $1\leq i\leq n-2$,}\\
&T_iT_j=T_jT_i,\quad\text{for $0\leq i<j-1\leq n-2$.}\end{aligned}
$$
Denote by $\Lambda$ the set of $m$-multipartitions of $n$. For
$\lam\in\Lambda$, let $M(\lam)$ be the set of standard
$\lam$-tableau. Then $\mathscr{H}_{n,m}$ has a cellular basis of the
form $\{m_{\fs \ft}^{\lam}\mid \lam\in\Lambda, \fs,\ft\in
M(\lam)\}$. See \cite{DG4} for details.

Let $L_{i}=q^{1-i}T_{i-1}\cdots T_{1}T_{0}T_{1}\cdots T_{i-1}$. Then
$L_{1}, L_{2}, \cdots, L_{n}$ is a family of Jucys-Murphy elements
of $\mathscr{H}_{n,m}$. If $i$ is in row $r$ column $c$ of the
$j$-th tableau of $\ft$, then $m_{\fs \ft}^{\lam}L_{i}\equiv
u_{j}q^{c-r}m_{\fs \ft}^{\lam}$. If $[1]_{q}\cdots
[n]_{q}\prod_{1\leq i<j\leq m}\prod_{|d|<n}(q^{d}u_{i}-u_{j})\neq 0$
and $q\neq 1$, then the Jucys-Murphy elements separate $M(\Lambda)$.
These were proved in \cite{JM}.

Denote $\mathscr{H}_{n,m}\otimes_{R}K$ by $\mathscr{H}_{n,m,K}$. The
following result has been proved in \cite{A} and \cite{RA}. We give
a new proof here.

\begin{thm}{\rm (\cite{A},\cite{RA})}
The center of $\mathscr{H}_{n,m,K}$ is equal to the set of symmetric
polynomials in the Jucys-Murphy elements if $[1]_{q}\cdots
[n]_{q}\prod_{1\leq i<j\leq m}\prod_{|d|<n}(q^{d}u_{i}-u_{j})\neq 0$
and $q\neq 1$.
\end{thm}
\noindent{Proof:}\, The algebra $\mathscr{H}_{n,m,K}$ satisfies the
conditions of the Proposition \ref{propjm} has been pointed out in
\cite{M}. By Theorem C, we only need to show that for any
$\lam,\mu\in\Lambda$ with $\lam\neq\mu$, $\{c_{\lam}(i)\mid 1\leq
i\leq n\}$ can not be obtained from $\{c_{\mu}(i)\mid 1\leq i\leq
n\}$ by permutations. Note that we can obtain these two sets by
$\ft^{\lam}$ and $\ft^{\mu}$ respectively.

{\it Case 1.} There exists $1\leq j\leq m$ such that $n_{j\lam}\neq
n_{j\mu}$. Then by the separation condition, the number of the
elements of the form $u_{j}q^{x}$ is $n_{j\lam}$ in
$\{c_{\lam}(i)\mid 1\leq i\leq M\}$ and is $n_{j\mu}$ in
$\{c_{\mu}(i)\mid 1\leq i\leq M\}$, where $x\in \mathbb{Z}$. This
implies that $\{c_{\lam}(i)\mid 1\leq i\leq M\}$ can not be obtained
from $\{c_{\mu}(i)\mid 1\leq i\leq n\}$ by permutations.

{\it Case 2.} $n_{j\lam}= n_{j\mu}$ for all $1\leq j\leq m$. Then
there must exist $1\leq s\leq m$, such that the partition of $n_{s}$
in $\lam$ is not equal to that in $\mu$ since $\lam\neq\mu$. Denote
the partitions by $\lam_{s}$ and $\mu_{s}$. Then the residues of
$\lam_{s}$ and $\mu_{s}$ are not the same. Now by Lemma \ref{lmp}
and the separation condition, the set of all the elements of the
form $u_{s}q^{x}$ in $\{c_{\lam}(i)\mid 1\leq i\leq M\}$ is
different from that in $\{c_{\mu}(i)\mid 1\leq i\leq n\}$. Then for
any $\lam,\mu\in\Lambda$ with $\lam\neq\mu$, $\{c_{\lam}(i)\mid
1\leq i\leq m\}$ can not be obtained from $\{c_{\mu}(i)\mid 1\leq
i\leq m\}$ by permutations.\qed

\noindent{\bf Remark.}\, If $q\neq 1$, then $[1]_{q}\cdots
[n]_{q}\prod_{1\leq i<j\leq m}\prod_{|d|<n}(q^{d}u_{i}-u_{j})\neq 0$
if and only if $\mathscr{H}_{n,m,K}$ is semisimple. See \cite{A} for
details.
\smallskip\bigskip

\bigskip\bigskip

\noindent{\bf Acknowledgments}\,\,\,The author acknowledges his supervisor Prof. C.C. Xi and the
support from the Research Fund of Doctor Program of Higher
Education, Ministry of Education of China. He also acknowledges Dr. Wei Hu and Zhankui Xiao for many helpful conversations.

\end{document}